\documentclass[11pt, oneside, reqno]{amsart}
\usepackage{fullpage}
\usepackage{amssymb,latexsym,amsmath,verbatim,layout,hyperref,amsthm,color}  
\numberwithin{equation}{section}
\usepackage[T1]{fontenc}
\usepackage[light,math]{anttor}
\usepackage {mathrsfs}                
  
\usepackage[colorinlistoftodos,prependcaption,textsize=tiny]{todonotes}   

\theoremstyle{mytheoremstyle}

\newtheorem{thm}{Theorem}[section]
\newtheorem{cor}{Corollary}
\newtheorem*{thm*}{Theorem}
\newtheorem{lem}[thm]{Lemma}
\newtheorem*{lem*}{Lemma}
\newtheorem{prop}[thm]{Proposition}
\newtheorem*{prop*}{Proposition}

\theoremstyle{definition}
\newtheorem{defn}[thm]{Definition}

\theoremstyle{remark}
\newtheorem*{rmk}{Remark}

\newtheorem{eg}{Example}
\usepackage{hyperref}

\renewcommand{\vec}[1]{\mathbf{#1}}

\newcommand{\sph}{\mathbb{S}}
\newcommand{\R}{\mathbb{R}}
\renewcommand{\P}{\mathbb{P}}
\newcommand{\E}{\mathbb{E}}
\newcommand{\N}{\mathbb{N}}

\title{Random Attractors for Stochastic Navier-Stokes equation on a 2D rotating sphere with stable L\'evy noise}
\author{Leanne Dong}
\date{
Faculty of Engineering and IT, Institute of Data Science, The University of Technology Sydney.\,\,	\today}
\begin{document}
\maketitle    
\begin{abstract}
In this paper we prove that the stochastic Navier-Stokes equations with stable L\'evy noise generates a random dynamical systems. Then we prove the existence of random attractor for the Navier-Stokes equations on 2D spheres under stable L\'evy noise (finite dimensional). We also deduce the existence of Feller Markov invariant measure.
\end{abstract}	


\section{Introduction and Motivation}
The study of asymptotic behaviour of dynamical systems has been one of the most fundamental problems in mathematical physics. One of the central notions is an attractor, which conveys crucial geometric information about the asymptotic regime of a dynamical system as $t\to\infty.$ It is well known that the 2D Navier-Stokes equations is dissipative and so have a global attractor, see for instance \cite{robinson2001infinite,temam1993infinite}. More precisely, there exists compact subset $K$ of the original phase space where all asymptotic dynamics fall in. Much of the theory of infinite dynamical systems devote to study the properties of this set $K$, which is called the \emph{global attractor} (see for instance \cite{temam1993infinite,robinson2001infinite}). For instance, one can show that, under certain mild conditions, it turns out that on $K$ one can define a group of solution operator $S(t)$ sensibly for all $t\in\R$, this defines a standard \emph{dynamical system},       
\begin{align*}    
    (K, \{S(t)_{t\in\R}\})
\end{align*}
A random (pullback) attractor is the pullback attractor where time-dependent forcing become random. Readers are referred to \cite{MR2844645} for a comparison of the three frameworks for the study of attractors, namely, attractors, pullback attractors and random attractors. Just like in the deterministic case, the theory of random attractor plays an important role in the study of the asymptotic behaviour of dissipative random dynamical system. The authors in \cite{MR1305587,MR1451294} developed a theory for the existence of random attractors for stochastic systems that closely comparable to the deterministic theory. Roughly speaking, a random attractor is a random invariant compact set which attracts every trajectory as time goes to infinity.
The strategy to prove the existence of random attractor is analogous to the method of proving global attractor in deterministic case, which involves two main methods. The first method requires one to find a bounded absorbing set and to prove asymptotic compactness of $S(t)$; The second method is to find compact absorbing set, and it turns out to be the method we employ in proving existence of random attractor in our case.
The existence theorems and continuous dependence of initial datum results in \cite{snse1} allows us to define a flow $\varphi:\R\times\Omega\times H\to H$ in the following sense:
\begin{align*}   \varphi_t(\omega)v_0=v(t,\omega;v_0),\quad t\in     T,\omega\in\Omega,\,\,v_0\in H.
\end{align*}  
In the paper \cite{BGQ17}, the authors proved the existence of random attractor for the Navier Stokes equations on 2D spheres under Gaussian-type forcing. To our awareness, there is only one result on these systems when they are perturbed by non-Gaussian noise   \cite{GU20141433}. Gaussian-type models has been used widely to model phenomenons in engineering , science and finance. The trajectories of a particle driven by Brownian motion are continuous in time almost surely, the displacement increases linearly in time (in mean square sense). Moreover, the probability density function decays exponentially in space. These models do not capture fluctuation with peculiar properties such as anomalous diffusion (See for instance \cite{BG90} and \cite{Y96}). 
A good candidate can be used to model complex phenomena involve irregular fluctuation with peculiar properties is L\'evy motion. Particularly non-Gaussian L\'evy motion, has been widely applied to Biology, Image processing, Climate forecast and certainly in Finance and Physics \cite{MR3616887,SZD17, ZCZ17,MR2356959,MR3457638}. Stable-type noise allows one to take into account at the same time noise with a large number of small random impulses and occasionally large random disturbance with infinite moment. From a fluid modeling point of view, although continuous models are good enough in a macroscopic scale, at an atomic scale, the model breaks down, and the use of a L\'evy process is compelling as fluid is not continuous at a microscopic scale  \cite{MMC11}. As a special non-Gaussian stochastic process, the stable-type process attracts more and more    mathematical interests due to the properties which the Gaussian process does not possess. The tail of Gaussian random variable decays exponentially which does not fit well for modeling processes with high variability or some extreme events, such as earthquakes or stock market crashes. In contrast, the stable L\'evy motion has a `heavy tail' that decays polynomially and can be useful for these applications. For instance, when heavier tails (relative to a Gaussian distribution) of asset returns are more pronounced,  the asymmetric $\alpha$-stable distribution becomes an appropriate alternative in modeling \cite{GSF12}.

The goal of this work is to investigate the dynamical behavior of the SNSE on 2D rotating spheres with additive stable L\'evy noise with $\beta\in (1,2)$.

\begin{align}\label{asnse1rds}
    du(t)=[-Au(t)-B(u(t),u(t))+\mathbf{C}u(t)+f]dt+\sum^m_{l=1}\sigma_l  dL_l(t)e_l,\quad u(0)=u_0,
\end{align}
where $A$ is the Stokes operator on the sphere defined as
    \begin{align}\label{optA}
         A : D(A)\subset H\to H, \quad  A =-P(\mathbf{\Delta}+2\text{Ric}),\quad D( A )=\mathbb{H}^2(\mathbb{S}^2)\cap V,
    \end{align}
The bilinear form $B: V\times V\to V'$ is defined by
\begin{align}\label{optB}    (B(u,v),w)=b(u,v,w)=\sum^3_{i,j=1}\int_{\sph^2}u_i\frac{\partial(v_k)_j}{\partial x_i}u_j dx,\quad w\in V.
\end{align}
With a slight abuse of notation, we denote $B(u)=B(u,u)$ and $B(u)=\pi(u,\nabla u)$.

The Coriolis operator $ \mathbf{C}_1 :\mathbb{L}^2(\mathbb{S}^2)\to\mathbb{L}^2(\mathbb{S}^2)$ is bounded linear on $\mathbb{L}^2(\mathbb{S}^2)$ and is defined by the formula\footnote{The angular velocity vector of earth is denoted as $\Omega$ in consistant to geophysical fluid dynamics Literature. It shall not be confused with the notation for probility space $\Omega$.}
\begin{align}
    ( \mathbf{C}_1 v)( x)=2\Omega( x\times  v( x))\text{cos}\theta,\quad  x\in\mathbb{S}^2.
\end{align} Moreover, we need a well defined bounded operator $\mathbf{C}=P\mathbf{C}_1$ in $H$. One can relate $\mathbf{C}$ and $\mathbf{C}_1$ as follows, 
\begin{align}\label{corio}
    ( \mathbf{C} u, u)=( \mathbf{C}_1  u,P u)=\int_{\mathbb{S}^2}2\Omega\text{cos}\theta(( x\times  u)\cdot u( x))dS( x)=0,\,\,\text{for}\,\,u\in H
\end{align}
Furthermore, $f\in H$ and $e_1,\cdots,e_m\in H$ are the eigenfunctions of the stokes operator $A$, $\{\sigma_l\}$ is a sequence of real numbers, $L_l(t)$, $(1\le l\le m)$ are mutually independent two-sided $\beta$-stable L\'evy processes $u=u(t,x,\omega)$ is now a random velocity of the fluid.
    
The goal in this work is in threefold.
\begin{itemize}
    \item Prove (\ref{asnse1rds}) generates a RDS $\varphi$;
    \item Establish the existence of random attractor for (\ref{asnse1rds});
    \item Establish the existence of a Feller Markov Invariant Measure supported by the random attractor.
\end{itemize}
To this end, we study the stationary ergodic solution of an Ornstein-Uhblenceck, make transformation to obtain some estimates of the solution respectively in space $H$ and $V$, then using the compact embedding of Sobolev space, we obtain the existence of compact random set which absorbs any bounded nonrandom subset of space $H$.

In section 2, we introduce some key terminology (namely L\'evy process, stable process, RDS, random attractors, Markov-invariant measures) to study random dynamics induced by our SNSE under jump noise. In section 3, we prove $\varphi$ indeed defines a random dynamical system along with a driving flow $\varphi$. This claim was accomplished by first identifying a suitable canonical sample probability space for (\ref{asnse1rds}) which ensures the linear stochastic Stokes equation remains stationary.
Then via an a priori estimate for a strong solution (bounded in $V$, compact in $H$) from our earlier work \cite{snse1}, we identified a compact absorbing set and consequently deduce the existence of a random attractor based on the assumption of a finite-dimensional noise. Finally, using the property of random attractor, we deduce the existence of random invariant measure which is supported by the random attractor of the spherical SNSE.
\section{Preliminary}\label{rdsdefns}
In this section, we first review some necessary notions and preliminary. The notion of L\'evy process follows closely with the book \cite{MR1739520}. The presentation here follows closely with \cite{MR2227218,MR1305587,arnold2013random} with some slight modification to Jump case based on the various papers on RDS under L\'evy noise \cite{MR2812588,Gu13}. 
The notion of random dynamical system is simply a generalisation of a deterministic dynamical system. In brief, an RDS has two features, one is the measurable dynamical system $\varphi$, which is used to model the random perturbations, and the other is the cocycle mapping $\vartheta$ defined over the dynamical system (see \cite{arnold2013random} for more detail). 

\subsection{Basic definitions}\label{basicdefn}
	\begin{defn}[L\'evy process]
An $H$-valued L\'evy process is a stochastic process $X=\{X(t),t\in [0,\infty)\}$ such that 
\begin{itemize}
	\item $X(0)=0$ a.s. and $X$ is stochastically continuous: $\forall\,\,\varepsilon>0$,
	\begin{align*}
		\lim_{t\downarrow 0}\P(|X(t)|>\varepsilon)=0.
	\end{align*}
	\item $X$ has independent increments, that is, $\forall\,\,0\le t_0<t_1<\cdots<t_n,$ the random vectors $X(t_0), X(t_1)-X(t_0), X(t_2)-X(t_1),\cdots, X(t_n)-X(t_{n-1})$ are independent,
	\item $X$ has stationary increments:
	\begin{align*}
		X_{t+s}-X_t\stackrel{d}=X_s\quad\forall\,s, t\ge 0,
	\end{align*}
	\item $t\mapsto X(t)$ is c\`adl\`ag a.s. .
\end{itemize}	
\end{defn}
	Note, that without the assumption of stationary increments, we have an \emph{additive process}.
We will need also a two-sided L\'evy process defined as follows. Let $X_1$ and $X_2$ be two independent L\'evy processes defined on the same probability space and with the same distribution. Then we define the two-sided L\'evy process 
\[X(t)=\left\{\begin{array}{lll}
X_1(t)&\mathrm{if}&t\ge 0\\
X_2(t)&\mathrm{if}&t<0
\end{array}\right.\quad t\in\R\,.\]
For a two-sided L\'evy process we will consider the filtration $\mathcal F_t=\sigma\left(X(s):\, s\le t\right)$ for all $t\in\R$.

The cylindrical L\'evy Noise used in this work is obtained by subordinating a cylindrical Wiener process by an arbitrary real valued, increasing L\'evy process. This increasing  L\'evy process is chosen to be  a $\frac{\beta}{2}$ stable (symmetric) process, with $\beta\in (0,2)$. Let us recall some basic facts from \cite{MR1280932}. First, recall that a real random variable $X$ is said to be $\beta$-stable with the, scale parameter $\sigma$, skewness parameter $\delta$, and shift parameter $\nu$, shortly $X\sim S_\beta\left(\sigma,\delta,\nu\right)$, 
if 
\[\E e^{i\theta X}=e^{i\theta\nu-|\sigma\theta|^\beta\left(1-i\delta c\,\mathrm{sgn}(\theta)\right)}\,,\]
where
\[c=\left\{\begin{array}{lll}
\left((\sigma\theta)^{1-\beta}-1\right)\tan\frac{\pi\beta}{2}&\mathrm{if}&\beta\neq 1\\
-\frac{2}{\pi}\log|\sigma\theta|&\mathrm{if}&\beta=1
\end{array}\right.\]
Note that in particular, $S_2(\sigma,0,\nu)=N(\nu,2\sigma^2)$ is Gaussian. \\
We have also the following definition. 
\begin{defn}
	A real valued random variable $X$ is said to be symmetric $\beta$-stable, $0<\beta\le 2$, if $X\sim S_\beta(\sigma,0,0)$ or, explicitly 
	\begin{align}\label{chf}
		\E e^{i\theta X}=e^{-\sigma^\beta|\theta|^\beta/2},\quad \theta\in\R.
	\end{align}
The name ``$\beta$-stable'' means that if $X_1,\cdots,X_m$ are independent and $\beta$-stable, then $\sum_{j\le m}\alpha_j X_j$ is $\beta$-stable, and
\begin{align*}
	\sigma (\sum_{j\le m}\alpha_j X_j)=
\left(\sum_{j\le m}|\alpha_j|^\beta\sigma(X_j)^p\right)^{1/\beta},
\end{align*}
which is obvious from (\ref{chf}).
\end{defn}
\begin{defn}
A random vector $X=(X_1,\cdots,X_N)$ with values in $\R^N$ is $\beta$-stable if each linear combination $\sum^N_{i=1}\alpha_i X_i$ is a real $\beta$-stable variable.\\
A random process $X=(X_t, t\in I)$ indexed by $I$ is called $\beta$-stable if for every $t_1,\cdots, t_N$ in $T$, $(X_{t_1},\cdots, X_{t_N})$ is a $\beta$-stable random vector. (p.131 in \cite{ledoux2013probability}, p.233 in \cite{talagrand2014upper})
\end{defn}
A natural generalisation of the $\R^n$ definition of stable L\'evy motion (see for instance p.113 \cite{MR1280932}) to the Hilbert space is the following
\begin{lem}\label{ben_lemma}
	A L\'evy process $\{X(t),t\ge 0\}$ on a Hilbert space is a $\beta$-stable L\'evy motion if and only if 
		$X(t)-X(s)\sim S_{\beta}((t-s)^{1/\beta},\delta,0)$
for some $0<\beta\le 2$, $-1\le\delta\le 1$.
\end{lem}
Now Let us recall the definition of random dynamical system (RDS) and cocylcle property.
\begin{defn}
A triple        $\mathfrak{T}=(\Omega,\mathcal{F},\vartheta)$ is said to be a \emph{measurable dynamical system} (DS) if $(\Omega,\mathcal{F})$ is a measurable space and $\vartheta:\R\times\Omega\ni (t,\omega)\mapsto\vartheta_t\omega\in\Omega$ is a $\mathcal B(\mathbb R)\otimes\mathcal F$-measurable map such that for all $t,s\in\R$, $\vartheta_{t+s}=\vartheta_t\circ\vartheta_s$.  A quadruple $\mathfrak{T}(\Omega,\mathcal{F},\P,\vartheta)$ is called a \emph{metric dynamical system} (RDS)  iff. $(\Omega,\mathcal{F},\P)$ is a probability space and $\mathfrak{T}':=(\Omega,\mathcal{F},\vartheta)$ is a measurable DS such that for each $t\in\R$, the map $\vartheta_t:\Omega\to\Omega$ is $\P$-preserving.
\end{defn}
\begin{defn}\label{defnrds}
    Given a metric DS $\mathfrak{T}$ and a Polish space $(X,d)$, a map $\varphi:\R\times\Omega\times X\ni (t,\omega,x)\mapsto\varphi(t,\omega)x\in X$ is called a \emph{measurable random dynamical system} (on $X$ over $\vartheta$), iff
\begin{itemize}
    \item $\varphi$ is $(\mathcal{B}(\R)\times\mathcal{F}\times\mathcal{B},\mathcal{B})$-measurable.
    \item The trajectories $\varphi(\cdot,\omega)x:\R\to X$ are c\`adl\`ag $\forall\,(\omega,x)\in\Omega\times\R$;
    \item $\varphi$ is $\vartheta$-cocycle:
\begin{align}\label{cocycle}
    \varphi(t+s,\omega)=\varphi(t,\vartheta_s\omega)\circ\varphi(s,\omega)\quad\forall\quad s,t\in\R,\,\,\varphi(0,\omega)=\text{id},\,\,\forall\,\,\omega\in\Omega.
\end{align}   
\end{itemize}


 It follows from the cocycle property that $\varphi(t,\omega)\cdot$ is automatically invertible. ($\forall t\in T$ and $\forall\,\P\,\,$ a.e. $\omega.$) In fact, $\varphi(t,\omega)^{-1}=\varphi(-t,\vartheta_t\omega)$ for $t\in T$. Instead of assuming (\ref{cocycle}) for all $\omega\in\Omega$ it suffices to assume (\ref{cocycle}) for all $\omega$ from a measurable $(\vartheta_t)_{t\in T}$-invariant subset of $\Omega$ of full measure.
\end{defn}
\subsection{Stochastic Calculus for two-sided time}
While we will assume our metric dynamical system has two sided time $T=\R$, in this subsection we briefly discuss the extension of stochastic calculus to two sided time. The material follows closely with section 2.3.2 \cite{arnold2013random}.

Let $(\Omega,\mathcal{F},P)$ denotes from now a complete probability space.
\begin{defn}[Two-Parameter Filtration, p.71 \cite{arnold2013random}] Assume $\mathcal{F}^t_s$, $s, t\in\R$, $s\le t$, is a two parameter family of sub $\sigma$-algebras of $\mathcal{F}$ with the following properties
\begin{itemize}
    \item $\mathcal{F}^t_s\subset\mathcal{F}^v_u$ for $u\le s\le t\le v$
    \item 
$\mathcal{F}^{t+}_s :=\cap_{u>t}\mathcal{F}^u_s=\mathcal{F}^t_s$, $\mathcal{F}^{t}_{s-} :=\cap_{u<s}\mathcal{F}^t_u=\mathcal{F}^t_s$ for $s\le t$,
    \item $\mathcal{F}^t_s$ contains all $\P$-null sets of $\mathcal{F}$ for every $s\le t$.
\end{itemize}
Then $\mathcal{F}^t_s$, $s\le t$ is called a (two-parameter) filtration on (on $(\Omega,\mathcal{F},P)$). We define
\begin{align*}
    \mathcal{F}^t_{-\infty}:=\vee_{s\le t}\mathcal{F}^t_s,\quad \mathcal{F}^{\infty}_{s}:=\vee_{t\ge s}\mathcal{F}^t_s.
\end{align*}
\end{defn}
\begin{defn}[Filtered DS, p.72 \cite{arnold2013random}]
    Let $(\Omega,\mathcal{F}^0,\P,\{\vartheta_t\}_{t\in\R})$ be a metric DS, let $\mathcal{F}$ be the $\P$-completion of $\mathcal{F}^0$, and let $\mathcal{F}^t_s$, $s\le t$, be a filtration in $(\Omega,\mathcal{F},\P)$. We call $(\Omega,\mathcal{F},\P,\{\vartheta\}_{t\in\R},\{\mathcal{F}^t_s\}_{s\le t})$ is filtered DS, if for all $s,t,u\in\R$, $s\le t$, we have
    \begin{align*}
        \vartheta^{-1}_u\mathcal{F}^t_s=\mathcal{F}^{t+u}_{s+u}.
    \end{align*}
\end{defn}
\subsection{Attraction and absorption}\label{ssec:attrabs}
For two random sets $A, B\subset X,$ we put
\begin{align*}
    d(A,B)=\sup_{x\in A}d(x,B)\quad\text{and}\quad\rho(A,B)=\max\{d(A,B),d(B,A)\}.
\end{align*} 
In fact, $\rho$  restricted to the family $\mathfrak{C}$ of all nonempty closed subsets on $X$ is a metric (see \cite{MR0467310}), and it is the so-called Hausdorff metric. From now on, let $\mathfrak{X}$ be the Borel $\sigma$-field on $\mathfrak{C}$ generated by open sets w.r.t. the metric $\rho$ \cite{MR1207308,MR0467310,  MR1993844}.
\begin{defn}
    Let us assume that $(\Omega,\mathcal{F})$ is a measurable space and let $(X,d)$ be a Polish space. A set-valued map $C:\Omega\to\mathfrak{C}(X)$ is said to be measurable iff. $C$ is $(\mathcal{F},\mathfrak{X})$-measurable. Such a map is often called a \emph{closed and bounded random set}. A closed and bounded random set $C$ will be called a \emph{compact random set} on $X$ if for each $\omega\in\Omega$, $C(\omega)$ is a compact subset of $X$.
\end{defn}
\begin{eg}
    A closed set valued map $K:\Omega\to 2^X$ is a random closed set.
\end{eg}
\begin{rmk}
    Let $f: X\mapsto\R$, be a continuous function on the Polish space $X$ and $\R:\Omega\mapsto\R$ an $\mathcal{F}$-measurable random variable. If the set $C_{f,R}(\omega) :=\{x:f(x)\le R(\omega)\}$ is nonempty for each $\omega\in\Omega$, then $C_{f,R}$ is a closed and bounded random set.
\end{rmk}
\begin{defn}
    Let $\varphi:\R\times\Omega\times X\ni(t,\omega,x)\mapsto\varphi(t,\omega)x\in X$ be measurable RDS on a Polish space $(X,d)$ over a metric $DS$ $\mathfrak{T}$. A closed random set $B$ is said to be $\varphi$ forward invariant iff. for all $\omega\in\Omega$,
    \begin{align*}
        \varphi(t,\omega)B(\omega)\subseteq B(\vartheta_t\omega)\quad\forall\quad t>0.
    \end{align*}
A closed random set $B$ is said to be \emph{strictly $\varphi$ invariant} iff. $\forall\,\omega\in\Omega$,
\begin{align*}
    \varphi(t,\omega)B(\omega)=B(\vartheta_t\omega)\quad\forall\quad t>0.
\end{align*}
\end{defn}
\begin{rmk}
    By substituting $\vartheta_{-t}\omega$ for $\omega$, we have the equivalent version of the above definition:
    \begin{align*}       \varphi(t,\varphi_{-t}\omega)B(\vartheta_{-t}\omega)\subseteq B(\omega),\quad\forall\,\,t>0,
    \end{align*}
    \begin{align*}       \varphi(t,\varphi_{-t}\omega)B(\vartheta_{-t}\omega)= B(\omega),\quad\forall\,\,t>0.
    \end{align*}
\end{rmk}

\begin{defn}\label{omegalimit}
    For a given closed random set $B$, the $\Omega$-limit set of $B$ is defined to be the set 
    \begin{align*}
        \Omega(B,\omega)=\Omega_{B}(\omega)=\bigcap_{T\ge 0}\overline{\bigcup_{t\ge T}\varphi(t,\varphi_{-t}\omega)B(\vartheta_{-t}\omega)}
    \end{align*}
\end{defn}
\begin{rmk}
\begin{itemize}
    \item [(i)] A priori $\Omega(B,\omega)$ can be an empty set.
    \item [(ii)] One has the following equivalent version of Definition \ref{omegalimit}:
\begin{align*}
    \Omega_B(\omega)=\{y: \exists t_n\to\infty, \{x_n\}\subset B(\vartheta_{-t_n}\omega), \lim_{n\to\infty}\varphi(t_n, \vartheta_{-t_n}\omega)x_n=y \}.
\end{align*} 
    \item [(iii)] Since $\overline{\bigcup_{t\ge T}\varphi(t,\vartheta_{-t}\omega)}$ is closed, $\Omega_B(\omega)$ is closed as well.
\end{itemize}
\end{rmk}

Given a probability space, a random attractor is a compact random set, invariant for the associated RDS and attracting every bounded random set in its basis of attraction. More precisely,
\begin{defn} A random set $A:\Omega\to\mathfrak{C} (X)$ is a \emph{random attractor} iff
\begin{itemize}
    \item $A$ is a compact random set;
    \item $A$ is $\varphi$-invariant, i.e. $\P$-a.s.
    \begin{align}
        \varphi(t,\omega)A(\omega)=A(\vartheta_t\omega),
    \end{align}
    \item $A$ is 
    attracting, in the sense that, for all $B\in X$ it holds
\begin{align*}
    \lim_{t\to\infty}\rho(\varphi(t,\vartheta_{-t}\omega)B(\varphi_{-t}\omega),A(\omega))=0.
\end{align*}

\end{itemize}
    
\end{defn}
The random attractor $A$ in the present work shall not be confused with the Stokes operator $A$ introduced in \cite{snse1}.

Let us now state a result on the existence of a random attractor, which is an extension of the Gaussian noise case in the pioneering work in \cite{MR1305587} to the general L\'evy noise case.

\begin{thm}\label{compactattract}
    Let $\varphi$ be a continuous in space, but c\`adl\`ag in time RDS on $X$. Assume the existence of a compact random set $K$ absorbing every deterministic bounded set $B\subseteq H$. Then there exists a random attractor $A$ given by
    \begin{align*}   
        A(\omega)=\overline{\bigcup_{B\subseteq X,      B\,\,\text{bounded}}\Omega_B(\omega)},\quad\omega\in\Omega.
    \end{align*}
\end{thm}
               
\begin{proof}   
    The proof is analogous to the proof of Theorem 3.11 in \cite{MR1305587}.
\end{proof}  
\subsection{Invariant measures on random sets}\label{ssec:imr}
In the final section, we prove the existence of invariant measure for the RDS $\varphi$ (Put in another way, the existence of random invariant measure). 
Let us now discuss the notion of \emph{random invariant measure}. 

\begin{defn}\label{imrds}
    Let $\varphi$ be a given RDS over a metric DS $\mathfrak{T}.$ A probability measure $\mu$ on $(\Omega\times X,\mathcal{F}\times\mathcal{B})$ is said to be an invariant measure for $\varphi$ iff.
\begin{itemize}
    \item $\Theta_t$ preserves $\mu: \Theta_t(\mu)=\mu$ for all $t>0$;
    \item The first marginal of $\mu$ is $\P$, i.e. $\pi_{\Omega}(\mu)=\P$ where 
    \begin{align*}
        \pi_{\Omega}: \Omega\times X\ni (\omega,x)\mapsto\omega\in\Omega.
    \end{align*}
\end{itemize}
\end{defn}
The following corollary gives the existence of invariant measure for a RDS $\varphi$. The proof follows from the Markov-Katutani fixed point theorem (see p.87 \cite{MR1993844} for more detail.).
\begin{cor}[p.374,\cite{MR1305587}]\label{cor44}
     Let $\varphi$ be an RDS, and suppose $\omega\mapsto A(\omega)$ is a compact measurable forward invariant set for $\varphi$. Then there exist invariant measure for $\varphi$ which are supported by $A$.
\end{cor}
Alternatively, one can construct the random invariant measure more explicitly via an Krylov-Bogoliubov type argument; we refer readers to p.87 \cite{MR1993844}.
\subsubsection{Markov Invariant Measures}

Based on the conditions in Theorem \ref{compactattract}, it is clear the attractor is measurable with respect to the past $\mathcal{F}^{-}$, since $\Omega_B$ is measurable for any nonrandom $B$.
  
Define two $\sigma$-algebra corresponding to the future and the past, respectively by
\begin{align*}
    \mathcal{F}^{+}=\sigma\{\omega\mapsto\varphi(\tau,\vartheta_t\omega):\tau,t\ge 0\},
\end{align*}
and
\begin{align*}
    \mathcal{F}^{-}=\sigma\{\omega\mapsto\varphi(\tau,\vartheta_{-t}\omega):\tau,0\le\tau\le t\}.
\end{align*}
Then $\vartheta^{-1}_t\mathcal{F}^+\subset\mathcal{F}^+$ for all $t\ge 0$ and $\vartheta^{-1}_t\mathcal{F}^-\subset\mathcal{F}^-$ for all $t\le 0$.

\begin{prop}
    Suppose $\omega\mapsto A(\omega)$ is an $\varphi$-invariant compact set which is measurable with respect to the past $\mathcal{F}^{-}$ for an RDS $\varphi$. Then there exist invariant measures $\mu$  supported by $A$ such that also $\omega\mapsto\mu_{\omega}$ is measurable with respect to $\mathcal{F}^{-}$.
\end{prop}
\begin{cor}[p.374\cite{MR1305587}]\label{markovmea}
     Under the conditions of the Proposition suppose in addition that $\varphi$ is an RDS whose one-point motions form a Markov family, and such that $\mathcal{F}^{+}$ and $\mathcal{F}^{-}$ are independent. Then there exists an invariant measure $\rho$ for the associated Markov semigroup. Furthermore, the limit 
    \begin{align*}
        \mu_{\omega}=\lim_{t\to\infty}\varphi(t,\vartheta_{-t}\omega)\rho
    \end{align*}
exists $\P$ a.s., $\rho=\int \mu_{\omega}d\P(\omega)=\E(\mu_{\cdot})$, and $\mu$ is a Markov measure.    
\end{cor}
\subsubsection{Feller Markov Invariant  measures}
By Corollary \ref{cor44} for an given RDS $\varphi$ on a Polish space $X$, one can find an invariant probability measure if an invariant compact random set $K(\omega)$, $\omega\in\Omega$ can be identified.
Hence Corollary \ref{cor44} is generalised as the following.
\begin{cor}
    A continuous in space, c\`adl\`ag in time RDS which has an invariant compact random set $K(\omega)$, $\omega\in\Omega$ has at least one invariant probability measure $\mu$ in the sense of definition \ref{imrds}.
\end{cor} 
One of the desirable property of solutions of stochastic PDE is Feller property. Let us now define a Feller invariant measure for a Markov RDS $\varphi$.  If $f:X\to\R$ is bounded Borel measurable function, then put 
\begin{align}\label{ptf}
    (P_t f)(x)=\E f(\varphi(t,x)),\quad t\ge 0, x\in X.
\end{align}
It is clear that $P_t f$ is also a bounded and borel measurable function. Moreover, one has the following result.
\begin{prop}\label{3.8}
    Assume that that RDS $\phi$ is a.s. contiuous in space for every $t\ge 0$. Then the  family $(P_t,t\ge 0)$ is Feller, i.e. $P_t f\in C_b(X)$ if $f\in C_b(X).$ Moreover, if the RDS $\varphi$ is c\`adl\`ag in time, then for any $f\in C_b(X)$, $(P_t f)(x)\to f(x)$ as $t\downarrow 0$.
\end{prop}
\begin{proof}
    For the first part, let us fix $t>0.$ If $x_n\to x$ in $X$, then it follows from the space continuity of $\varphi(t,\cdot)$ that $(P_t f)(x_n)\to (P_t f)(x)$ using the Lebesgue dominated convergence theorem.

For the second part, note that  for a given $x\in X$ from the c\`adl\`ag property of $\varphi(\cdot,x,\omega): [0,\infty)\to X$ for a.e. $\omega$ it follows that one has $(P_t f)(x)\to f(x)$ as $t\to 0$ if $x\in X.$ 
\end{proof}
A RDS $\varphi$ is called Markov iff the family $(P_t, t\ge 0)$ is a semigroup on $C_b(X)$, that is, $P_{t+s}=P_t\circ P_s$ for all $t, s\ge 0.$

\begin{defn}
    A Borel probability measure $\mu$ in $H$ is said to be invariant w.r.t. $P_t$ if 
    \begin{align*}
    P^*_t\mu:=\int_X P_t(x,\Gamma)\mu(dx)=\mu(\Gamma),\quad\forall\Gamma\in\mathcal{B}(X),\quad \forall t
    \end{align*}
where $(P^*_t)(\Gamma)=\int_H P_t(x,\Gamma)\mu(dx)$ for $\Gamma\in\mathcal{B}(H)$ and $P_t(x,\cdot)$ is the transition probability, $P_t(x,\Gamma)=P_t(1_{\Gamma})(x)$
\end{defn}
\vspace{1cm}
Finally, a Feller invariant probability measure for a Markov RDS $\varphi$ on $H$ is, by definition, an invariant probability measure for the semigroup $(P_t, t\ge 0)$ define by (\ref{ptf}).

In view of Corollary \ref{markovmea},  if a Markov RDS $\varphi$ on a Polish space $H$ has an invariant compact random set $K(\omega), \omega\in\Omega$, then there exists a Feller invariant probability measure $\mu$ for $\varphi$. More precisely we have the following corollary.
   
\begin{cor}\label{fpm}
    If a c\`adl\`ag time and space continuous RDS $\varphi$ has an invariant compact random set, then there exists a feller invariant probability measure $\mu$ for $\varphi.$
\end{cor}
     
\section{Random Dynamical systems generated by the SNSE on a rotating unit sphere}\label{randomdynamic}
Having established the well-posedness in the earlier work \cite{snse1}, we are in a position to define an RDS $\varphi$ corresponding to the problem 
\begin{align}\label{asnse4}
	du(t)+Au(t)+B(u(t),u(t))+\mathbf{C}u=fdt+GdL(t),\quad u(0)=u_0,
\end{align}
where $L$ is an $H$-valued stable L\'evy process and $G:H\to H$ is a bounded operator.

But first, we need to determine a sample (canonical) probability space for which the dynamics of the driving noise remains stationary.

\subsection{Some analytic facts}\label{ssec:someanalytic}
Recall that $X=\mathbb{L}^4(\sph^2)\cap H$ denote the Banach space endowed with the norm
\begin{align*}
    |x|_X=|x|_H+|x|_{\mathbb{L}^4(\sph^2)}.
\end{align*}

Recall Assumption 1 in \cite{snse1}, namely,
the space $K\subset H\cap \mathbb{L}^4$ is a Hilbert space such that for any $\delta\in (0,1/2)$,
\begin{align}
    A^{-\delta}: K\to H\cap \mathbb{L}^4(\sph^2)\quad\text{is $\gamma$-radonifying. }
\end{align}
This assumption is satisfied if $K=D(A^s)$ for some $s>0.$
\begin{rmk}
    Under the above assumption the space $K$ can be taken as the RKHS of the cylindrical Wiener process  $W(t)$ on $H\cap \mathbb{L}^4$.
\end{rmk}
\vspace{1cm}
 Let $(\Omega,\mathcal{F},\P)$ be a complete probability space, where $\Omega=D(\R,E)$ of c\`adl\`ag functions defined on $\R$ take value in $E$ with the following Skorohod metric
\begin{align*}
    d(l_1,l_2)&=\sum^{\infty}_{i=1}(1\wedge d^{\circ}_i(l_1,l_2))\quad\forall\,\,l_1,l_2\in D, 
\end{align*}
where $l^i_1(t):=g_i(t)l_1(t)$ and 
$l^i_2(t):=g_i(t)l_2(t)$ with
\begin{align*}
    g_i(t)&:=
\begin{cases}
    1, &\text{ if }|t|\le i-1\\
    i-t,&\text{ if }i-1\le |t|\le i\\
    0,&\text{ if }|t|\le i
\end{cases}
\end{align*}
\begin{align*}   d^{\circ}_i(l^i_1,l^i_2)=\inf_{\lambda\in\Lambda}\left(\sup_{-i\le s<t\le i}|\log\frac{\lambda(t)-\lambda(s)}{t-s}|\vee\sup_{-i\le t\le i}|l_1(t)-l_2(\lambda(t))|\right),
\end{align*}
where $\Lambda$ denotes the set of strictly increasing, continuous function $\lambda(t)$ from $\R$ to itself with $\lambda(0)=0$. This skorohod space is a complete separable metric space which is taken as the canonical sample space. Let $\mathcal{F}$ be the Borel $\sigma$-algebra of the Polish space $(D(\R,X),d)$. For every $t\in\R$ we have the evaluation map $L_t: D(\R,X)\to\R$ denote by $L_t(\omega)=\omega(t)$. Then we have $\mathcal{F}=\sigma(L_t, t\in\R)$, that is, $\mathcal{F}$ is the smallest  $\sigma$-algebra generated by the family of maps $\{L_t: t\in\R\}$. 
Let $\P$ be the unique probability measure which makes the canonical process a two-sided L\'evy process with paths in $D(\R;E)$, that is, $\omega(t)=L_t(\omega)$.

Define the shift
\begin{align*}    (\vartheta_t\omega)(\cdot)=\omega(t+\cdot)-\omega(t)\quad t\in\R,\omega\in\Omega.
\end{align*}
Then the map $(t,\omega)\to\vartheta_t(\omega)$ is continuous and measurable \cite{arnold2013random} and the (L\'evy) probability measure $\P$ is $\vartheta$ invariant, that is, $\P(\vartheta^{-1}_t(T))=\P(T)$ for all $T\in\mathcal{F}$. This flow is an ergodic dynamical system with respect to $\P$. Thus $(\Omega,\mathcal{F},\P,(\vartheta)_{t\in\R})$ is a metric DS.

\subsection{Ornstein-Uhlenbeck process}\label{ssec:ourds}
In the following subsections we are concerned with the linear evolutionary Stokes equations. 

Recall, the equation
\begin{align*}
    \begin{cases}
        \dot{u}(t)+Au(t)=f(t), &\text{  }t\in[0,T],\\
        u(0)=u_0.
    \end{cases}
\end{align*}
If $A$ generates a $C_0$-semigroup in a Banach space $E$ and $f:[0,T]\to E$ is such a function that 
\[\int_0^T|f(t)|_E\,dt<\infty,\]
then the solution is given by 
\begin{align*}
    u(t)=e^{-tA}u_0+\int^t_0 e^{-\alpha(t-s)A}f(s)ds.
\end{align*}
In particular, we have
\begin{prop}\label{transform}
Let $L$ be a L\'evy process taking value in $E$, such that for any $T>0$
\[\int_0^T|L(t)|_E\,dt<\infty\,.\]
 Then the solution of the differential equation
    \begin{align*}
        \dot{V}(t)+\alpha V(t)=L(t),\quad V(0)=0,\quad \alpha>0
    \end{align*}
is given by 
\begin{align*}
    V(t)=\int^t_0 e^{-\alpha(t-s)A}L(s)ds.
\end{align*}
\end{prop}

\subsubsection{Stochastic convolution and integrating by parts}
Here we quote a useful integration by part formula from \cite{Xu:2012aa} which allows us to attain the desired regularity for which the RDS $\varphi$ exist. 

Consider the following Ornstein-Uhlenbeck process generated by the Stokes operator on $\sph^2$,
\begin{align*}
    z_t=\int^t_0 e^{-A(t-s)}GdL(s)=\sum^{\infty}_{l=1}z_l e_l,
\end{align*}
where $\{e_l: l=1,\cdots\}$ is the complete orthonormal system of eigenfunctions of $A$ in $H$ and 
\begin{align}
    z_t(t)=\int^t_0 e^{-\lambda_l(t-s)}\sigma_l dL^l(s),
\end{align} 
where $\lambda_l$ are the eigenvalues of the Stokes operator $A$.
By the It\^o product formula, see Theorem 4.4.13 of  \cite{applebaum2007levy} for any $l\ge 1$, one has that
\begin{align*}
    L^l(t)=\int^t_0\lambda_l e^{-\lambda_l(t-s)}L^l(s)ds+\int^t_0e^{-\lambda_l(t-s)}dL^l(s)+\int^t_0 \lambda_l e^{-\lambda_l(t-s)}\Delta L^l(s)ds,
\end{align*}
where $\Delta L^l(s)=L^l(s)-L^l(s-)$. Since $L^l(t)$ is a $\beta$ stable process, $\Delta L^l(s)=0$ a.e. for $s\in[0,t]$ and so we have
\begin{align*}
    \int^t_0\lambda_l e^{-\lambda_l(t-s)}\Delta L^l(s)ds=0.
\end{align*}
Therefore,
\begin{align*}
    z_l(t)=\sigma_l L^l(t)-\int^t_0\lambda_l e^{-\lambda_l(t-s)}\sigma_l L^l(s)ds.
\end{align*}
Hence, if we assume that $\sigma_l=0$ for $l>m$ for a certain finite $m>1$ then
\begin{align*}
    z(t)=L(t)-Y(t),
\end{align*}
where 
\begin{align*}
    Y(t)=\int^t_0 A e^{-A(t-s)}L(s)ds.
\end{align*}
In this case we clearly have
\begin{align*}
    Y(t)\in H\quad\text{a.s.}\quad (t\ge 0).
\end{align*}
\subsubsection{Regularity of Shifting flow}
To prove our stochastic Navier-Stokes system generates a RDS, we will transform it into a random PDE in $X$ with aid of the integration by part technique introduced earlier. We need to give a meaning to the Ornstein-Uhlenbeck process in the metric dynamical system $(\Omega,\mathcal{F},\P,\{\vartheta_t\}_{t\in\R})$ given by 
\begin{align*}
z(\vartheta_t\omega):=\hat{z}(t)=\int^t_{-\infty}\widehat A^{1+\delta}
e^{-(t-r)\widehat A}(\tilde{\omega}(t)-\tilde{\omega}(r))dr,\quad t\in\R.
\end{align*}
Our goal now is to show $\widehat{z}(t)$ is a well defined element in $X:=\mathbb{L}^4(\sph^2)\cap H$ for a.e. $\omega$.
But first, we need the following couple of results, which can be viewed as generalisation of Theorem 4.1 and Theorem 4.4 in \cite{MR2584982} to the case where the Ornstein Uhlenbeck  generator is  $\hat{A}+\alpha I$ in place of $A$, where
\begin{align*}
    \widehat{A}&=\nu A+C,\quad D(\hat{A})=D(A),\quad A=-P(\vec{\Delta}+2\text{Ric}).
\end{align*}
Recall that 
\begin{align}\label{eq_anal1}
    |\widehat{A}^{1+\delta}e^{-t(\hat{A}+\alpha I)}|_{\mathcal{L}(X,X)}\le C t^{-1-\delta}e^{-(\mu+\alpha) t},\quad t > 0.
\end{align}
\begin{prop}\label{flowthm} Assume $\beta\in (1,2)$, $p\in(0,\beta)$ and 
\[\sum_{l=1}^\infty\left|\sigma_l\right|^\beta\lambda_l^{\beta/2}<\infty\,.\]
Then 
\begin{equation}\label{eq1}   \E\int^t_{-\infty}|\hat{A}e^{-(t-r)\hat{A}}(\tilde{\omega}(t)-\tilde{\omega}(r))|_X^p dr<\infty.
\end{equation}
Moreover, for $\P$ \emph{almost every} $\tilde{\omega}\in D(\R,X)$, $t\in\R$ the function 
\begin{align}\label{eq_conv1}  \hat{z}(t)=\hat{z}(\tilde{\omega})(t)=\int^t_{-\infty}\hat{A}e^{-(t-r)\hat{A}}(\tilde{\omega}(t)-\tilde{\omega}(r))dr,\quad t\in\R
\end{align}
is well defined and c\'adl\'ag in $X$. 
Furthermore, for any $\kappa>0$ such that $\kappa p>1$ there 
 exists a random variable $C$ depending on $\beta$, $p$, $\sigma$, $\delta$ such that 
\begin{align}\label{eq_power}   |\hat{z}(\tilde{\omega})(t)|_X\le C(\beta, p,\sigma,\delta,\tilde{\omega})(1+|t|^{\kappa}).
\end{align}
\end{prop}

\begin{proof}
\emph{Part I}\hspace{1cm} We will show first that the L\'evy process $L$ is is c\`adl\`ag in $X$ and 
\begin{equation}\label{omega}
\mathbb E\sup_{t\le T}|L(s)|^p_X<\infty\,.
\end{equation}
By Lemma 3.4 in \cite{snse1}, the process $L$ is c\'adl\'ag in $V$, hence in $\mathbb L^2\left(\mathbb S^2\right)$ and finally in $X$. It remains to show \eqref{omega}. Recall Lemma 3.1, that is,
    \begin{align}
        \E\sup_{t\le T}|A^{\delta}L(t)|^p\le C(\beta,p)\left(\sum_{l\ge 1}|\sigma_l|^{\beta}\lambda^{\beta\delta}_l\right)^{\frac{p}{\beta}} t^{\frac{p}{\beta}}<\infty.
    \end{align}
Putting $\delta=0$, we have,
\begin{align}
    \E\sup_{t\le T}|L(t)|^p\le C(\beta,p)\left(\sum_{l\ge 1}|\sigma_l|^{\beta}\right)^{\frac{p}{\beta}} t^{\frac{p}{\beta}}<\infty,
\end{align}    
and putting $\delta=\frac{1}{2}$ we have,
    \begin{align}
        \E\sup_{t\le T}|A^{\frac{1}{2}}L(t)|^p\le C(\beta,p)\left(\sum_{l\ge 1}|\sigma_l|^{\beta}\lambda^{\beta/2}_l\right)^{\frac{p}{\beta}} t^{\frac{p}{\beta}}<\infty.
    \end{align}
So
\begin{align*}
    \E\sup_{t\le T}|L(s)|^p_X&\le c\E\sup_{t\le T}|L(s)|^p+c\E\sup_{t\le T}|L(s)|^p_{\mathbb{L}^4(\sph^2)}\\
    &\le c\E\sup_{t\le T}|L(s)|^p+C\E\sup_{t\le T}\left(|L(s)|^{\frac{p}{2}}|L(s)|^{\frac{p}{2}}_V\right)\quad\text{via Ladyzhenskaya inequality}\\
    &\le c\E\sup_{t\le T}|L(s)|^p+C\E\sup_{t\le T}|L(s)|^p_V\quad\text{via Poinc\'are inequality in \cite{snse1}}\\
    &=C(\beta,p)\left(\sum_{l\ge 1}|\sigma_l|^{\beta}\right)^{p/\beta} s^{\frac{p}{\beta}}+C(\beta,p)\left(\sum_{l\ge 1}|\sigma_l|^{\beta}\lambda^{\beta/2}_l\right)^{\frac{p}{\beta}} s^{\frac{p}{\beta}}\\
    &\le C(\beta,p)\left(\sum_{l\ge 1}|\sigma_l|^{\beta}\lambda^{\beta/2}_l\right)^{\frac{p}{\beta}} s^{\frac{p}{\beta}}.
\end{align*}
\emph{Part II }\hspace{1cm} In what follows we use the fact that $\tilde\omega(t)=L(t)$ $\mathbb P$-a.s. Using the change of variables  $s=t-r$, we obtain 
\begin{align*}
    \E\int^t_{-\infty}|\hat{A}e^{-(t-r)\hat{A}}(\tilde{\omega}(t)-\tilde{\omega}(r))|^p dr&=\int^{\infty}_0\E|\hat{A}e^{-s\hat{A}}(\tilde{\omega}(t)-\tilde{\omega}(t-s))|_X^p ds\\
    &=\int^{\infty}_0\E|\hat{A}e^{-s\hat{A}}\tilde{\omega}(s)|_X^p ds.\\
\end{align*}
Using \eqref{eq_anal1} with $\gamma=\alpha+\mu$ we have  
\begin{align*}
        \int^{\infty}_0\E|\hat{A}e^{-s\hat{A}}\tilde{\omega}(s)|_X^p ds&\le C\int^{\infty}_0 \frac{e^{-p\gamma s}}{s^{p}}\E|A^\delta\tilde{\omega}(s)|_X^p ds\\
    \le &\
C    \int^{\infty}_0\frac{e^{-p\gamma s}}{s^{p}}C(\beta,p)s^{\frac{p}{\beta}}(\sum_{l\ge 1}|\sigma_l|^{\beta}\lambda^{\beta\delta}_l)^{\frac{p}{\beta}}ds
<\infty,
\end{align*}
since $p-\frac{p}{\beta}<1$ and we infer that $\widehat z(t)$ is well defined in $X$ $\mathbb P$-a.s. using the same arguments as in the proof of \eqref{omega} above. \\
We will prove \eqref{eq_power}. Applying  Lemma 6.6 in \cite{snse1}, with the Banach space $B=X$ and $\kappa$ such that $\kappa p>1$ we obtain 
\begin{align}\label{eq_power1}
    |\widehat{z}(\tilde{\omega})(t)|_X\le C(\beta, p,\sigma,\delta,\kappa,\tilde{\omega})(1+|t|^{\kappa}),
\end{align}
and \eqref{eq_power} follows. \\
\emph{Part III}\quad 
One has to check $\hat{z}$ is right continuous with left limit in $X$. To this end note first that 
\[\begin{aligned} \widehat z(t)&=\int_{-\infty}^t \widehat Ae^{-(t-s)\hat{A}}\left(\omega(t)-\omega(s)\right) ds\\
&=\left(\widehat A \int_0^\infty e^{-s\hat{A}} ds\right)\omega(t)-\int_{-\infty}^t \widehat Ae^{-(t-s)\hat{A}}\omega(s)\, ds\\
&=\omega(t)-\int_{-\infty}^t \widehat Ae^{-(t-s)\hat{A}}\omega(s)\, ds,
\end{aligned}\]
since $\widehat A$ is invertible. The function $\omega$ is c\`adl\`ag in $X$ by assumption. We will show that the function 
\[F(t,\omega)=\int_{-\infty}^t \widehat Ae^{-(t-s)\hat{A}}\omega(s)\, ds\]
is continuous in $X$ for $\P$-a.e. $\omega$. Indeed, for $s,t\in\mathbb R$ such $r<t$ we have 
\[\begin{aligned}
\int_{-\infty}^t \widehat Ae^{-(t-s)\widehat{A}}\omega(s)\, ds&=\int_{-\infty}^r \widehat Ae^{-(t-s)\widehat{A}}\omega(s)\, ds+\int_{r}^t \widehat Ae^{-(t-s)\widehat{A}}\omega(s)\, ds\\
&=\widehat Ae^{-(t-r)\widehat{A}}\int_{-\infty}^r e^{-(r-s)\widehat{A}}\omega(s)\, ds+\int_{r}^t \widehat Ae^{-(t-s)\widehat{A}}\omega(s)\, ds\\
&=\widehat Ae^{-(t-r)\widehat{A}}h+I(t)\,.
\end{aligned}\]
Since the semigroup $e^{-s\widehat{A}}$ is analytic, we find that the function $t\to\widehat Ae^{-(t-r)\widehat{A}}h$ is continuous for $t>r$. Let us consider $I(t)$. By Sobolev embeddings we have a continuous embedding $\mathbb H^{1,2}\subset\mathbb W^{\frac14,4}$. Therefore for $\delta$ small enough the function $t\to A^\delta\omega(t)$ is locally bounded in $\mathbb L^4$ for a.e. $\omega$. Then 
\[I(t)=\int_r^t \widehat A^{1-\delta}e^{-(t-s)\widehat{A}}\widehat A^\delta\omega(s)\, ds\]
is continuous for $t>r$, 
 again by standard properties of analytic semigroups. 

\end{proof}

\begin{thm}\label{shift}
    Under the assumption of Proposition \ref{flowthm}, for  $\mathbb P$-a.e. $\omega\in D(\R,X)$, 
    \begin{align*}
        \hat{z}(\vartheta_s\omega)(t)=\hat{z}(\omega)(t+s),\quad t, s\in \R.
    \end{align*}
\end{thm}
\begin{proof}
    The proofs of the first three parts follows from closely from Theorem 4.8 and Proposition 8.4 in \cite{MR2584982}, see also Theorem 9 in \cite{huang2013random}. For the last part, since $(\vartheta_s \omega)(r)=\omega(r+s)-\omega(s)$, $r\in\R$, we have
    \begin{align*}
        \hat{z}(\vartheta_s\omega)(t)&=\int^t_{-\infty}Ae^{-(t-r)A}[\vartheta_s\omega(t)-\vartheta_s\omega(r)]dr\\
        &=\int^t_{-\infty}\widehat Ae^{-(t-r)A}[\omega(t+s)-\omega(r+s)]dr\\               
        &=\int^{t+s}_{-\infty}\widehat Ae^{-(t+s-r')A}[\omega(t+s)-\omega(r')]dr'=\hat{z}(\omega)(t+s).
    \end{align*}
\end{proof}
Now, put $(\tau_s\zeta)(t)=\zeta(t+s)$, $t,s\in\R.$ Therefore $\tau_s$ is linear, bounded map from $D(\R,X)$ into $D(\R,X)$. Moreover, the family $(\tau_s)_{s\in\R}$ is a $C_0$ group on $D(\R,X).$ 
Hence the shifting property could be re-expressed as
\begin{cor} For $\mathbb P$-a.e. $\omega\in D(\mathbb R,X)$ 
    For $s\in\R$, $\tau_s\circ\hat{z}=\hat{z}\circ\vartheta_s$, that is
    \begin{align*}  \tau_s(\hat{z}(\omega))=\hat{z}(\vartheta_s(\omega)),\quad\omega\in D(\R;X).
    \end{align*}
\end{cor}

\begin{prop}
    The process 
    \begin{align*}   z_{\alpha}(t)=\int^t_{-\infty}e^{-(t-s)(\hat{A}+\alpha I)}dL(s),
\end{align*}
where the integral is 
intepreted in the sense of \cite{MR2584982} is well defined and is identified as a solution to the equation 
    \begin{align*}
  dz_{\alpha}(t)+(\hat{A}+\alpha I)z_{\alpha}dt=dL(t),\quad t\in\mathbb R.
    \end{align*}
The process $z_{\alpha}$, $t\in\R$ is a stationary OU process. 
\end{prop}

We define
\begin{align*}   z_{\alpha}(\omega):=\hat{z}(\hat{A}+\alpha I;\omega)\in D(\R,X),
\end{align*}
i.e. for any $t\ge 0$,
\begin{align*}
    z_{\alpha}(\omega)(t) :=\int^t_{-\infty}(\hat{A}+\alpha I)e^{-(t-r)(\hat{A}+\alpha I)}(\omega(t)-\omega(s))ds
\end{align*}
By Proposition \ref{transform},
\begin{align*}    \frac{d^+}{dt}z_{\alpha}(\omega)(t)+(\hat{A}+\alpha I)\int^t_{-\infty}(\hat{A}+\alpha I)^{1+\delta}e^{-(t-r)(\hat{A}+\alpha I)}(\omega(t)-\omega(s))ds=L(t).
\end{align*}
Therefore $z_{\alpha}(t)$ satisfies 
\begin{align}\label{za}
        \frac{d^+}{dt}z_{\alpha}(t)=(\hat{A}+\alpha I)z_{\alpha}=\dot{\omega}(t),t\in\R.
\end{align}
It follows from Theorem \ref{shift} that
\begin{align*}
    z_{\alpha}(\vartheta_s\omega)(t)=z_{\alpha}(\omega)(t+s),\quad\omega\in D(\R,X),\,\,t,s\in\R.
\end{align*}
Similar to our definition of L\'evy process $L_t$, i.e. $L_t(\omega):=\omega(t)$, we can view the ODE as a definition of $z_{\alpha}(t)$ on $(\Omega,\mathcal{F},\P)$, equation (\ref{za}) suggests that this process is an Ornstein Uhblenck process.
               
Now we have enough tools to prove the cocycle property of RDS, and 
this allows us to prove $(\varphi,\vartheta)$ is an RDS. The proof follows same lines as  Theorem 6.15 in \cite{MR2227218}.
\subsection{Random dynamical system generated by the SNSE on a sphere with L\'evy noise}\label{rdsconstr}

Let us fix $\alpha\ge 0$ and put $\Omega=\Omega(E)$.

We define a map $\varphi=\varphi_{\alpha}:\R\times\Omega\times H\to H$ by
\begin{align*}
    (t,\omega,x)\mapsto v(t,\hat{z}_{\alpha}(\omega))(x-\hat{z}_{\alpha}(\omega)(0))+\hat{z}_{\alpha}(\omega)(t).
\end{align*}
In what follows, write $z=z_{\alpha}$ for simplicity.

Put in another way,
\begin{align*}
    \varphi=\varphi_{\alpha}(t,\omega)x&:=v(t,z_{\alpha}(\omega))(x-z_{\alpha}(\omega)(0))+z_{\alpha}(\omega)(t)\\
    &=u(t,\omega;x)\quad\forall\,\,t\in T,\omega\in\Omega,\,x\in H,
\end{align*}
where $u(\cdot;\omega,u_0)$ is the solution of the integral equation corresponding to given $\omega\in\Omega$, $u_0\in H$ and $\varphi$ satisfies the definition of RDS.

Since $\varphi(t)=\varphi(t,\vartheta_t(\omega))v_0$ and $v(0)=v_0$. Then $\varphi(0,\omega)=I.$ It is clear that $\varphi(0,\omega)=I.$
Because $z(\omega)\in D(\R; X)$, $z(\omega)(0)$ is a well-defined element of $H$ and hence $\varphi$ is well defined. Furthermore, we have the first main result of this work.
\begin{thm}\label{t5rds}
    $(\varphi,\vartheta)$ is a random dynamical system.
\end{thm}
To prove the claim, one simply check the definition of a random dynamical system (see subsection \ref{defnrds}). 

\begin{proof}
First, we check \emph{Measurability}. Suppose $u_0\in V$ and $t\in T$ is fixed, the map $\omega\mapsto\varphi(t,\omega)u_0\in H$ is measurable because the solution $u(t,\omega;u_0)$ is constructed as the (poinwise in $\omega$) limit of sucessive approximation of the contraction, which is measurable being explicitly defined in term of measurable objects. Finally, if $u_0\in H$, then $u(t,\omega;u_0)$ is the limit of $u(t,\omega;u^0_n)$ with $u^0_n\in V$. The required measurability is assured.

Next, we check \emph{Continuous dependence on initial data}.

The proof follows the similar line of the proof of uniqueness for the strong solution in \cite{snse1}.

Third, we check 
\emph{C\`adl\`ag property of $\varphi(t,\omega)$}. This turns out to be an easy task: The c\`adl\`ag property of $\varphi_t u_0$ is clear from the proof of existence and uniqueness of the solution $u$.

Lastly, we will check the \emph{cocycle property} of $\varphi$, namely, for any $x\in H$, one has to check,
 \begin{align}\label{cocyc} \varphi(t+s,\omega)x=\varphi(t,\vartheta_s\omega)\varphi(s,\omega),\quad t,s\in\R.
    \end{align}
From the definition of $\varphi$, noting from the cocycle property of $z$,  $\hat{z}(\vartheta_s\omega)(t)=\hat{z}(\omega)(t+s)$,  $z(\omega)(s)=z(\vartheta_s\omega)(0)$ for all $s\in\R$, we have for all $t,s\in\R$,
\begin{align*}   \varphi(t+s,\omega)x=v(t+s,z(\omega)(t+s))(x-z(\omega)(0))+z(\omega)(t+s),
\end{align*}

\begin{align*}   \varphi(t,\vartheta_s\omega)\varphi(s,\omega)x&=v(t,z(\vartheta_s\omega)(t))(x-z(\omega)(0))+z(\omega)(s)-z(\vartheta_s\omega)(0)+z(\vartheta_s\omega)(t)\\
    &=v(t,z(\vartheta_s\omega)(t))(v(s,z(\omega))(s))(x-z(\omega)(0))+z(\vartheta_s\omega)(t).
\end{align*}
In view of (\ref{shift}), to prove (\ref{cocyc}), we need to prove
\begin{align*}
    v(t+s,z(\omega)(t+s))(x-z(\omega)(0))=v(t,z(\vartheta_s\omega)(t))(v(s,z(\omega)(0))).
\end{align*}
Now, fix $s\in\R$, define $v_1$, $v_2$ by
\begin{align*}
    v_1(t)&=v(t+s,z(\omega)(t+s))(x-z(\omega)(0))\quad t\in\R,\\
    v_2(t)&=v(t,z(\vartheta_s\omega)(t))(v(s,z(\omega)(s))(x-z(\omega)(0))),\quad t\in\R.
\end{align*}
Because $v(0,z(\vartheta_s\omega)(0))(x-z(\vartheta_s\omega)(0))=x-z(\vartheta_s\omega)(0)$, one infer that
\begin{align*}  v_1(0)&=v(s,z(\omega)(s))(x-z(\omega)(0))\\  &=v(0,z(\vartheta_s\omega)(0))(v(s,z(\omega)(s))(x-z(\omega)(0)))=v_2(0).
\end{align*}    
Since $\R\ni t\mapsto v(t,z(\omega))$ is a solution to 
\begin{align}\label{ODE}
    \begin{cases}
        \frac{dv}{dt+}=-\nu Av-B(v)-B(v,z)-B(z,v)-B(z)+\alpha z+f, \\
        v(0)=v_0.
    \end{cases}
\end{align}
On the other hand, in view of our earlier existence uniqueness results, the fact $v$ takes value in $D(A)$ implies that $v(t)$ is differentiable for almost every $t$. 
We have
\begin{align*}
    v'(t)&=-\nu Av_1(t+s,z(\omega)(t+s))-B(v_1(t+s,z(\omega)(t+s))+z(\omega)(t+s)+\alpha z(\omega)(t+s)+f\notag\\
    &=-\nu A v_1(t,z(\omega))-B(v_1(t,z(\omega))+z(\omega)(t+s))+\alpha z(\omega)(t+s)+f.
\end{align*}
On the other hand for $v_2$,
\begin{align*}
    \frac{dv(t,z(\vartheta_s\omega)(t))}{dt+}=-\nu Av(t,z(\vartheta_s\omega)(t))-B(v(t,z(\vartheta_s\omega)(t)))+z(\vartheta_s\omega)(t)+\alpha z(\vartheta_s\omega)(t)+f.
\end{align*}
Therefore, $v_1$, $v_2$ solve respectively
\begin{align*}
    \begin{cases}
        v'_1(t)=-\nu Av'_1-B(v'_1(t)+z(\omega)(t+s))+\alpha z(\omega)(t+s)+f,\\
        v_1(0)=v(z(\omega))(s)(x-z(\omega)(0)),
    \end{cases}
\end{align*}
\begin{align*}
    \begin{cases}
        v'_2(t)=-\nu Av'_2-B(v'_2(t)+z(\vartheta_s\omega)(t))+\alpha z(\vartheta_s\omega)(t)+f,\\
        v_1(0)=v(z(\omega))(s)(x-z(\omega)(0)).
    \end{cases}
\end{align*}
By cocycle property of $z$, $z(\vartheta_s\omega)(t)=z(\omega)(t+s)$ for $t\in\R$.
\end{proof}
Therefore, $v_1$, $v_2$ are solutions to (\ref{ODE}) with the same initial data $v(s,z(\omega)(s))(x-z(\omega)(0))$ at $t=0$. Then it follows from the uniqueness of solution to (\ref{ODE}) that $v_1=v_2$, $t\in\R.$

\subsection{Existence of random attractors}\label{rattractim}
This subsection aims to establish the existence of random attractors. The main lines follow from classical lines of proving global attractors by finding compact absorbing sets. However, as pointed out in the paper \cite{MR1305587}, the analysis of Navier-Stokes equations with additive noise in our case requires some non-trivial consideration. In particular, a critical question arised when analyzing the estimate $\frac{d^+}{dt}|v(t)|^2$,  the usual estimates for the nonlinear term $b(v(t),z(t),v(t))$ yields a term $|v(t)|^2 |z(t)|^4_4,$ so we were not able to deduce any bound in $H$ for $|v(t)|^2$ under the classical lines (see for instance section 6 in \cite{MR3306386}). Nevertheless, in light of the method developed in \cite{MR1305587}, via the usual change of variable and by writing the noise and the Ornstein-Uhlenbeck process as an infinite sequence of 1D processes, we are able to show there exist random attractors to our system \ref{asnse1rds} as well. In what proceed we will detail our proof. First we need a few Lemmas from our two companion papers.

\begin{lem}\label{asymlevy}[Lemma 3.2, \cite{snse1}]
	Suppose that there exists some $\delta>0$ such that $\sum_{l\ge 1}|\sigma_l|^{\beta}\lambda^{\beta\delta}_l<\infty$. Then for all $p\in (0,\beta)$,
	\begin{align}\label{hypo}
		\E|A^{\delta}L(t)|^p\le C(\beta,p)\left(\sum_{l\ge 1}|\sigma_l|^{\beta}\lambda^{\beta\delta}_l\right)^{\frac{p}{\beta}} t^{\frac{p}{\beta}}<\infty.
	\end{align}
\end{lem}

\begin{lem}\label{dom_A}[Lemma 1.7,\cite{snse3}]
We have 
\[\sup_{-1\le t\le 0}|Az(t)|^2<\infty\,.\]
\end{lem}

Now, using Lemma \ref{dom_A} and Lemma \ref{asymlevy} applied with $\delta=\frac12$ we find that the process $z$ is c\`ad\`ag in $V$ and   
\begin{align}\label{zzvaz}
    \sup_{-1\le t\le 0}\left(|z(t)|^2+|z(t)|^2_V+|Az(t)|^2\right)<\infty\quad\P\,\,\text{a.s.}\quad .
\end{align}
Using  equation (4.12) in \cite{MR2773026}, one can now choose $\alpha>0$ such that    
\begin{align}\label{z10}
    4\eta m\E|z_1(0)|\le \frac{\lambda_1}{4},
\end{align}
where $\lambda_1$ is the first eigenvalue of $A$, since $\E|z_1(0)|^p\to 0$ as $\alpha\to\infty.$

 From (\ref{z10}) and the Ergodic Theorem we obtain
\begin{align*} \lim_{t_0\to-\infty}\frac{1}{-1-t_0}\int^{-1}_{t_0}4\eta\sum^m_{l=1}|z_l(s)|ds=4\eta m\E|z_1(0)|<\frac{\lambda_1}{4}.
\end{align*}
Put $\gamma(t)=-\frac{\lambda_1}{2}+4\eta\sum^m_{l=1}|z_l(t)|$, we get
\begin{align}\label{gp1}
    \lim_{t_0\to-\infty}\frac{1}{-1-t_0}\int^{-1}_{t_0}\gamma(s)ds<-\frac{\lambda_1}{4}.
\end{align}
From this fact and by stationarity of $z_l$ we finally obtain
\begin{align}\label{gp2}   \lim_{t_0\to-\infty}e^{\int^{-1}_{t_0}\gamma(s)ds}=0\quad\P-\text{a.s.}\quad,
\end{align}

\begin{align}\label{gp3}   \sup_{t_0<-1}e^{\int^{-1}_{t_0}\gamma(s)ds}|z(t_0)|^2<\infty,\quad\P-\text{a.s.}\quad .
\end{align}

\begin{align}\label{gp4}
    \int^{-1}_{-\infty}e^{\int^{-1}_{\tau}\gamma(s)ds}(1+|z_l(\tau)|^2+|z_l(\tau)|^2_V+|z_l(\tau)|^2|z_l(\tau)|)d\tau<\infty,\quad\P-\text{a.s.}\quad. 
\end{align}
for all $1\le j$, $l\le m.$ Indeed, note for instance that for $t<0$,
\begin{align*} \frac{z_l(t)}{t}=\frac{z_l(0)}{t}-\frac{1}{t}(\alpha+A_l)\int^0_t z_l(s)ds+\frac{L_l(t)}{t},
\end{align*}
whence $\lim_{t\to-\infty}\frac{z_l(t)}{t}=0$ $\P$-a.s., which implies (\ref{gp2}) and (\ref{gp3}).
Consider the abstract SNSE
\begin{align*}   du+[Au+B(u)+\mathbf{C}u]dt=fdt+GdL(t)
\end{align*}
and the Ornstein-Uhlenback equation
\begin{align*}
    dz+(\hat{A}+\alpha I) z dt=GdL(t),
\end{align*}
where $L(t)=\sum^{m}_{l=1}e_l L_l(t)$.
We now use the change of variable $v(t)=u(t)-z(t).$ Then, by subtracting the Ornstein-Uhlenback equation from the abstract SNSE, we find that $v$ satisfies the equation
\begin{align}\label{odeim}
    \frac{d^+v}{dt}=-\nu Av(t)-\mathbf{C}v(t)-B(u,u)+f+\alpha z.
\end{align}

Recall the Poincare inequalities
\begin{align}\label{poincarerds}
    |u|^2_V&\ge\lambda_1|u|^2,\quad\forall\quad u\in V,\\
    |Au|^2&\ge\lambda_1|u|^2,\quad\forall\quad u\in D(A).
\end{align}

\begin{lem}\label{ineqattractor}
    Suppose that $v$ is a solution to problem (\ref{ODE}) on the time interval $[t_0,\infty)$ with $z\in L^4_{\text{loc}}(\R,\mathbb{L}^4(\mathbb{S}^2))\cap L^2_{\text{loc}}(\R,V')$ and $ t_0\ge 0.$ Then, for any $t\ge\tau\ge t_0$, one has
\begin{align}\label{ineqrds}
    |v(t)|^2\le |v(\tau)|^2 e^{\int^t_{\tau}\gamma(s)ds}+\int^t_{t_0}e^{\int^t_s\gamma(\xi)d\xi}2p(s)ds,
\end{align}
where 
\begin{align}\label{pt} p(t)=c|f|^2+c\alpha|z|^2+\delta|z|^2\sum^m_{l=1}|z_l(t)|,
\end{align}
\begin{align}\label{gt}
    \gamma(t)=-\frac{\lambda_1}{2}+4\delta\sum^m_{l=1}|z_l(t)|
\end{align}
for all $t_0\le \tau\le t$
and $c$ depends only on $\lambda_1$
\end{lem}
\begin{proof}
    The proof will be provided shortly.
\end{proof}

Let $H$, $A: D(A)\subset H\to H$, $V=D(A^{1/2})=D(\hat{A}^{1/2})$ and $B(u,v):V\times V\to V', \,\,Cu$ be spaces and operators introduced in the previous section. Suppose that there exists a constant $c_B>0$ such that 
\begin{align}\label{buvw0rds}
    \langle B(u,v),w\rangle=|b(u,v,w)|&\le c_B|u|^{1/2}|u|^{1/2}_V|v|^{1/2}|v|^{1/2}_V|w|_V,\quad\forall\quad u,v,z\in V,
\end{align}

\begin{align*}
    \langle B(u,v),v\rangle\le c_B|u|^{1/2}|Au|^{1/2}|v|_V|z|
\end{align*}
for all $u\in D(A)$, $v\in V$ and $z\in H.$
 Moreover, let $f\in H$, $e_1,\cdots,e_m\in H$ be given, $\{\sigma_l\}$ is a sequence of real numbers. Consider \ref{asnse1rds} again,
\begin{align*}
    du(t)=[-Au(t)-B(u(t),u(t))+Cu(t)+f]dt+\sum^m_{l=1}\sigma_l  dL_l(t)e_l,\quad u(0)=u_0.
\end{align*} 
As in \cite{snse1}, assume that $e_l$ are the eigenfunctions of the stoke operator $A$, $1\le l\le m$, there exists $\delta>0$ such that 
\begin{align}\label{bsumrds}
    |\langle B(u,e_l),u\rangle|\le\delta|u|^2,\quad u\in V, l=1,\cdots,m.
\end{align}
\begin{rmk}
    In bounded domain or in $\mathbb{S}^2$, one has
    \begin{align}\label{assumpB}
            \langle B(u,e_l),u\rangle=\sum^3_{i,j=1}\int_{\mathbb{S}^2}u_i\frac{\partial(e_l)_j}{\partial x_i}u_j dx.
    \end{align}
\end{rmk}
In this case assumption (\ref{bsumrds}) is satisfied when $e_l$ are Lipschitz continuous in $\mathbb{S}^2.$ Put $L(t)=\sum^m_{l=1}e_l L_l(t).$

\subsubsection{Stochastic flow}
Consider the abstract SNSE
\begin{align*}
    du+[Au+B(u)+Cu]dt=fdt+GdL(t),
\end{align*}
and the Ornstein-Uhlenback equation
\begin{align*}
    dz+(\hat{A}+\alpha I) z dt=GdL(t).
\end{align*}
From the discussion from the earlier subsubsection, it is clear that $z(t)$ is a stationary ergodic solution with continuous trajectories take value in $D(A)$. So we can transform the SNSE to a random PDE. The main advantage is that we can solve the equation $\omega$-wise due to the absence of the stochastic integral.
  
We now use the change of variable $v(t)=u(t)-z(t).$ Then, by subtracting the Ornstein-Uhlenback equation from the abstract SNSE, we find that $v$ satisfies the following equation
\begin{align}\label{odevrds}
    \frac{dv}{dt+}=-\nu Av(t)-Cv(t)-B(u,u)+f+\alpha z.
\end{align}
    
Now recall the following theorem from \cite{snse1}, namely,
\begin{thm}\label{t45}
    Assume that $\varphi_l\in D(A)$, $1\le l\le m$, there exists $\eta>0$ such that     
\begin{align}\label{bsum}
    |\langle B(u,\varphi_l),u\rangle|\le\eta|u|^2,\quad u\in V, l=1,\cdots,m
\end{align} is satisfied. Then for $\P$-a.s. $\omega\in\Omega$, there hold
\begin{itemize}
    \item For all $t_0\in\R$ and all $v_0\in H$, there exists a unique solution $v\in C([t_0,+\infty];H)\cap L^2_{\text{loc}}([t_0,+\infty);V)$ of equation (\ref{odevrds}) with initial value $v_0$.
    \item If $v_0\in V$, then the solution belongs to $C([t_0,+\infty);V)\cap L^2_{\text{loc}}([t_0,+\infty);D(A))$.
    \item hence, for every $\varepsilon>0, v(t)\in C([t_0+\varepsilon,+\infty);V)\cap L^2_{\text{loc}}([t_0+\varepsilon,+\infty);D(A))$.
    \item Denoting the solution by $v(t,t_0;\omega,v_0)$, then the map $v_0\mapsto v(t,t_0;\omega,v_0)$ is continuous for all $t\ge t_0,$ $v_0\in H$.
\end{itemize}
\end{thm}
Now Let us define the transition semigroup for the flow $\varphi$ as 
\begin{align*}
    P_t f(x)=\E f(\varphi(t,x)).
\end{align*}
\begin{cor}
    It follow from Theorem \ref{t45} the transition semigroup for the Markov RDS $\varphi$ has Feller property in $H$. That is, $P_t: C_b(H)\to C_b(H)$
\end{cor}

    Having the map $v_0\mapsto v(t,t_0;\omega,v_0)$, where $v(t,t_0;\omega,v_0)$ is the solution to (\ref{odevrds}) with $v(t_0)=v_0,$ we can now define a stochastic flow $\varphi(t,\omega)$ in $H$ by setting 
    \begin{align*}
        \varphi(t,\omega)u_0=v(t,0;\omega,u_0-z_{\alpha}(\omega)(0))+z_{\alpha}(\omega)(t).
    \end{align*}

\subsubsection{Absorbing in $H$ at time $t=-1$ }
In what proceed, assume $\omega\in\Omega$ is fixed; the results will hold $\P$ a.s.. Suppose $t_0\le -1$ and $u_0\in H$ be given, and let $v$ be the solution of equation (\ref{odevrds}) for $t\ge t_0$, with $v(t_0)=u_0-z(t_0,\omega)$ (which was denoted above by $v(t,0;\omega,u_0-z(0,\omega))$). Using the well known identity $\frac{1}{2}\partial_t|v(t)|^2=(v(t),v(t))$, and the assumption $\langle B(u,v),v\rangle =0$ and the antisymmetric term $(Cv,v) =0$ we have
\begin{align}\label{deq5}
    \frac{1}{2}\frac{d^+}{dt}|v|^2&=-\nu(Av,v)-\langle B(u,z),u\rangle+(\alpha z,v)+\langle f,v\rangle\\
    &\le -\nu |v|^2_V-\langle B(u,z),u\rangle+\alpha |z||v|+|f||v|.
\end{align}
By the definition of $z$ and assumptions (\ref{bsum}),
\begin{align*}
\left|\langle B(u,z),u\rangle\right|&=\left|\sum^m_{l=1}\langle B(u,e_l),u\rangle e_l\right|\le\delta|u|^2\sum^m_{l=1}|z_l|\\
    &\le 2\delta|v|^2\sum^m_{l=1}|z_l|+2\delta|z|^2\sum^m_{l=1}|z_l|.
\end{align*}
and the inequalities
\begin{align*}
    ( \alpha z, v ) = c\alpha |z|^2+c'|v|^2,
\end{align*}

\begin{align*}
    \langle f, v\rangle \le c |f|^2+c'|v|^2.
\end{align*} 
For simplicity we take $\nu=1$.   
Then via Young inequality, one can show that there exists $c, c'>0$ depending only on $\lambda_1$ such that 

\begin{align*}
    \frac{1}{2}\frac{d^+}{dt}|v|^2+\frac{1}{2}|v|^2_V\le \left(-\frac{\lambda_1}{4}+2\delta\sum^m_{l=1}|z_l|+2c'\right)|v|^2+c|f|^2+c\alpha|z|^2+2c|z|^2_V+2\delta|z|^2\sum^m_{l=1}|z_l|.
\end{align*}
Let $\gamma(t)$, and $p(t)$ are defined as in \ref{ineqattractor}. Namely,

\begin{align*}
    p(t)=c|f|^2+c\alpha|z|^2+\delta|z|^2\sum^m_{k=1}|z_k(t)|,
\end{align*}
\begin{align*}
    \gamma(t)=-\frac{\lambda_1}{2}+4\delta\sum^m_{l=1}|z_l(s)|,
\end{align*}
 we have 
 \begin{align}\label{dineqrds}     \frac{1}{2}\frac{d^+}{dt}|v|^2+\frac{1}{2}|v|^2_V\le \frac{1}{2}\gamma(t)|v|^2+p(t),
 \end{align}
\begin{align*}
    \frac{d^+}{dt}|v(t)|^2\le\gamma(t)|v(t)|^2+2p(t).
\end{align*}
Invoking Gronwall Lemma over the interval $[a,\infty)$, we have (\ref{ineqrds}).
\begin{lem}
    There exists a random radius $r_1(\omega)>0$ such that for all $\rho>0$ there exists (a deterministic) $\bar{t}\le -1$ such that the following holds $\P$-a.s. For all $t_0\le \bar{t}$ and for all $u_0\in H$ with $|u_0|\le\rho$, the solution $v(t,t_0;\omega,u_0-z(s))$ of equation (\ref{ODE}) over $[t_0,\infty]$ with $v(t_0)=u_0-z_{\alpha}(t_0)$ satisfies the inequality
    \begin{align*}
        |v(-1,t_0;\omega,u_0-z_{\alpha}(t_0,\omega))|^2\le r^2_1(\omega).
    \end{align*}
\end{lem}
\begin{proof}
Apply Lemma \ref{ineqattractor} with $t=-1$, $\tau=t_0$, we have

\begin{align}\label{esth}
    |v(-1)|^2&\le |v(t_0)|^2 e^{\int^{-1}_{t_0}\gamma(\xi)d\xi}+\int^{-1}_{t_0}e^{\int^t_{s}\gamma(\xi)d\xi}2p(s)ds\notag\\   
    &\le 2e^{\int^{-1}_{t_0}\gamma(\xi)d\xi}|u_0|^2+2e^{\int^{-1}_{t_0}\gamma(\xi)d\xi}|z(t_0)|^2+\int^{-1}_{-\infty}e^{\int^t_{s}\gamma(\xi)d\xi}2p(s)ds.
\end{align}
    Put 
    \begin{align}
        r^2_1(\omega)&=2+2\sup_{t_0\le -1}e^{\int^{-1}_{t_0}\gamma(\xi)d\xi}|z(t_0)|^2+\int^{-1}_{-\infty}e^{\int^t_{s}\gamma(\xi)d\xi}2p(s)ds
    \end{align}
which is finite $\P$ a.s. due to the stationarity of $z_l$ (namely, equations (1.20) and (1.21) in \cite{snse3}).
   
So, given $\rho>0$, choose $\bar{t}$ such that 
\begin{align*}         
    e^{\int^{-1}_{t_0}\gamma(\xi)d\xi}\rho^2\le 1
\end{align*}
for all $t_0\le\bar{t}$. The claim then follows from (\ref{ineqrds}). We remark that $t_0$ depending on $\omega.$
\end{proof}  

Taking $t\in[-1,0]$ and $\tau=-1$ in (\ref{ineqrds}) we have
\begin{align*}
        |v(t)|^2\le |v(-1)|^2 e^{\int^t_{-1}\gamma(\xi)d\xi}+\int^t_{-1}e^{\int^t_s\gamma(\xi)d\xi}2p(s)ds.
\end{align*}
Let us now come back to (\ref{dineqrds}): 
\begin{align*}
 \frac{d^+}{dt}|v|^2+|v|^2_V\le \gamma(t)|v|^2+2p(t).
\end{align*}

Integrate over $[-1,0]$,
\begin{align*}
    |v(0)|^2-|v(-1)|^2+\int^0_{-1}|v(s)|^2_Vds\le \left(\int^0_{-1}\gamma(\xi)d\xi\right)\left(\sup_{-1\le t\le 0}|v(t)|^2\right)+\int^0_{-1}2p(s)ds.
\end{align*}
Therefore,
\begin{align*}
    \int^0_{-1}|v(s)|^2_Vds\le |v(-1)|^2+ \left(\int^0_{-1}\gamma(\xi)d\xi\right)\left(\sup_{-1\le t\le 0}|v(t)|^2\right)+\int^0_{-1}2p(s)ds.
\end{align*}
Therefore, from above lemma we deduce
\begin{lem}
    There exists two random variables $c_1(\omega)$ and $c_2(\omega)$ depending on $\lambda_1$, $e_1,\cdots,e_m$ and $|f|$ such that for all $\rho>0$ there exists $\bar{t}(\omega)\le -1$ such that the following holds $\P$ a.s. $\forall\,t_0\le\bar{t}$ and for all $u_0\in H$ with $|u_0|\le\rho$, the solution $v(t,\omega;t_0,u_0-z(t_0,\omega))$ of equation (\ref{odevrds}) over $[t_0,\infty]$ with $v(t_0)=u_0-z(t_0)$ satisfies
    \begin{align*}
        |v(t,\omega;t_0,u_0-z(t_0,\omega))|^2\le c_1(\omega)\quad\forall\quad t\in[-1,0],
    \end{align*}
    
    \begin{align*}       \int^0_{-1}|v(s,\omega;t_0,u_0-z(t_0,\omega))|^2_Vds\le c_2(\omega).
    \end{align*}
\end{lem}
\begin{proof}
    Put 
    \begin{align*}
        c_1(\omega)=e^{\int^t_{-1}\gamma(\xi)d\xi}r^2_1(\omega)+\int^t_{-1}e^{\int^t_{s}\gamma(\xi)d\xi}p(s)ds,
    \end{align*}
    
    \begin{align*}
        c_2(\omega)=r^2_1(\omega)\left(1+\int^0_{-1}\gamma(\xi)d\xi\right)+\int^0_{-1}2p(s)ds,
    \end{align*}
with $r_1(\omega)$ as in the previous lemma. Then, given $\rho>0$, it suffices to choose $t(\omega)$ as in the proof of that previous lemma.
\end{proof}
\subsubsection{Absorption in $V$ at $t=0$}
 From (\ref{ODE}) we have (by multiplying $Av$ left and right and noting $(v_t,Av)=\frac{1}{2}\frac{d^+}{dt}|v|^2_V$), using inequality (\ref{buvw0rds}), and use the Young inequality with $ab=\sqrt{\frac{1}{e}}a\sqrt{e}b$, $p=2$ 
 \begin{align*}
     ab\le\frac{a^2}{2}+\frac{b^2 }{2},
 \end{align*}
 With the choice of $e=\frac{\nu}{4}$, one has that
 \begin{align*}    
     \langle f, Av\rangle\le \frac{2}{\nu}|f|^2+\frac{\nu}{8}|Av|^2,
 \end{align*}
 
 \begin{align*}
     \langle \alpha z, Av\rangle\le \frac{2}{\nu}|\alpha z|^2+\frac{\nu}{8}|Av|^2,
 \end{align*}

 \begin{align*}     c_B|u|^{1/2}|Au|^{1/2}|u|_V|Av|\le 2\nu c^2_B|u||Au||u|^2_V+\frac{\nu}{8}|Av|^2,
 \end{align*}

 \begin{align*}     \frac{1}{2}\frac{d^+}{dt}|v|^2_V&=-\nu|Av|^2-
\langle B(u,u),Av\rangle +\langle f, Av\rangle+\langle \alpha z, Av\rangle    \\
&\le -\nu|Av|^2+c_B|u|^{1/2}|Au|^{1/2}|u|_V|Av|+|f||Av|+|\alpha z||Av|\\
&\le -\frac{5\nu}{8}|Av|^2+2\nu c^2_B|u||Au||u|^2_V+\frac{2}{\nu}(|f|^2+|\alpha z|^2)
\intertext{With $q(t)=\frac{2}{\nu}(|f|^2+|\alpha z|^2)$ and noticing $|Au|\le |Av|+|Az|$, }
&\le -\frac{5\nu}{8}|Av|^2+2\nu c^2_B|u||Av||u|^2_V+2\nu c^2_B|u||Az||u|^2_V+q(t)\\
\intertext{Apply Young inequality with $e=\nu/2$ for the nonlinear term,}
&\le -\frac{5\nu}{8}|Av|^2+\frac{\nu}{8}|Av|^2+4\nu c^4_B|u|^2|u|^2_V|u|^2_V+2c^2_B|u||Az||u|^2_V+q(t)\\
&\le -\frac{\nu}{2}|Av|^2+4\nu c^4_B|u|^2|u|^2_V|u|^2_V+2c^2_B|u||Az||u|^2_V+q(t)\\
&\le -\frac{\nu}{2}|Av|^2+8\nu c^4_B|u|^2|u|^2_V|v|^2_V+8\nu c^4_B|u|^2|u|^2_V|z|^2_V+2\nu c^2_B|u||Az||u|^2_V+q(t).
\end{align*}   
 Temporarily disregard the $|Av|^2$ term, we have
 \begin{align*}
     \frac{d^+}{dt}|v|^2_V\le 16\nu c^4_B|u|^2|u|^2_V|v|^2_V+16\nu c^4_B|u|^2|u|^2_V|z|^2_V+4\nu c^2_B|u||Az||u|^2_V+2q(t).
 \end{align*}
Invoking Gronwall Lemma, we get for $s\in[-1,0]$,
\begin{align*}
    |v(0)|^2_V&\le e^{\int^0_s 16\nu c^4_B|u(\tau)|^2|u(\tau)|^2_Vds}\times\left(|v(s)|^2+\int^0_s(16\nu c^4_B|u|^2|u|^2_V|z|^2_V+4\nu c^4_B|u||Az||u|^2_V+2q(t))d\sigma\right).
\end{align*}
Integrating in $s$ over $[-1,0]$ we obtain
\begin{align}\label{v0V}
    |v(0)|^2_V&\le \left(\int^0_{-1}|v(\tau)|^2d\tau+\int^0_{-1}[16\nu c^4_B|u|^2|u|^2_V|z|^2_V+4\nu c^4_B|u||Az||u|^2_V+2q(t)]d\sigma\right)e^{\int^0_{-1}16\nu c^4_B|u(\tau)|^2|u(\tau)|^2_Vd\tau}.
\end{align}

\begin{lem}
    There exists a random radius $r_2(\omega)>0$, depending only on $\lambda_1$, $e_1,\cdots,e_m$ and $f$,   such that for all $\rho>0$ there exists (a deterministic) $\bar{t}\le -1$ such that the following holds $\P$-a.s. For all $t_0\le \bar{t}$ and for all $u_0\in H$ with $|u_0|\le\rho$, the solution $v(t,t_0,\omega;u_0-z(t_0,\omega))$ of equation \ref{odevrds} over $[t_0,\infty]$ with $v(t_0)=u_0-z_{\alpha}(t_0)$, put $u(t,t_0,\omega;u_0)=z(t,\omega)+v(t,t_0,\omega;u_0-z(t_0,\omega))$. Then
    \begin{align*}
        |u(0,\omega;t_0,u_0-z_{\alpha}(t_0,\omega))|^2_V\le r^2_2(\omega).
    \end{align*}
\end{lem}
\begin{proof}
In view of (\ref{v0V}), we need to estimate the term.
\begin{align*}
    \int^0_{-1}|u(s)|^2|u(s)|^2_V ds.
\end{align*}
Now using the fact $u=v+z$ and $|u|^2\le 2|v|^2+2|z|^2=2(|v|^2+|z|^2)$,, then the two terms can be estimated as following.
\begin{align*}
&\    \int^0_{-1}|u(s)|^2|u(s)|^2_V ds\\
&\le \sup_{-1\le t\le 0}|u(t)|^2\int^0_{-1}|u|^2_V ds\\
&\le \sup_{-1\le t\le 0} 2(|v(t)|^2+|z(t)|^2)\left(\int^0_{-1}2|v(\tau)|^2+2|z(\tau)|^2\right)\\
&\le 2(c_1(\omega)+\sup_{-1\le t\le 0}|z(t)|^2)2(c_2(\omega)+\int^0_{-1}|z(s)|^2 ds)\\
&= 2c_3(\omega)2c_4(\omega),
\end{align*}    

\begin{align*}
&\    \int^0_{-1}|u(s)||u(s)|^2_V ds\\
&\le \sup_{-1\le t\le 0}|u(t)|\int^0_{-1}|u(s)|^2_Vds\\
&\le (c_1(\omega)+\sup_{-1\le t\le 0}|z(t)|)2c_4(\omega).
\end{align*}
Hence, put 
\begin{align*}
    c_3(\omega)&=c_1(\omega)+\sup_{-1\le t\le 0}|z(t)|^2,\\
    c_4(\omega)&=c_2(\omega)+\int^0_{-1}|z(t)|^2_Vds,\\
    c_5(\omega)&=c_1(\omega)^{1/2}+\sup_{-1\le t\le 0}|z(t)|.
\end{align*}
Then, (\ref{v0V}) becomes
\begin{align*}
    |u(0)|^2_V&\le 2|z(0)|^2_V+2|v(0)|^2_V\\
&\le 2|z(0)|^2_V\\&+2\left[c_2(\omega)+64 c^4_B c_3(\omega)c_4(\omega)\sup_{-1\le t\le 0}|z(t)|+8c^2_B    c_5(\omega)c_4(\omega)\sup_{-1\le t\le 0}|Az(t)|+\int^0_{-1}2q(s)ds\right] e^{64 c^4_B c_3(\omega)c_4(\omega)} \\
&=: r^2_2(\omega).
\end{align*}
\end{proof}
Hence there exists a random ball in $V$ which absorbs the bounded sets of $H.$ Since $V$ is compactly embedded in $H$, there exists a compact set $K\subset H$ such that, for all bounded set $B\subset H$ there exists $\bar{t}\le -1$ such that $\varphi B\subset K$ $\P$ almost surely.

\subsection{Existence of Feller Markov Invariant Measures}\label{mim}
In this subsection, we prove the existence of random attractor implies the existence Feller Markov invariant measures.

\begin{thm}\label{t5attr}
    The stochastic flow associated with the SNSE with additive L\'evy noise (\ref{asnse1rds}) has a compact random attractor, in the sense of Theorem \ref{compactattract}. Moreover, the Markov semigroup induced by the flow of $H$ has an invariant measure $\rho$ in the sense of Corollary \ref{markovmea}. The associated flow-invariant Markov measure $\mu$ on $H\times\Omega$ has the property that its disintegration $\omega\mapsto\mu_{\omega}$ is supported by the attractor.
\end{thm}
\begin{proof}
    Recall that, in the language of the stochastic flow associated with our SNSE (\ref{asnse1rds}),
    \begin{align*}
    u(0,\omega;t_0,u_0)=\varphi(t_n,\vartheta_{-t_0}\omega)u_0=v(0,\omega;t_0,u_0-z(s))+z(t).
    \end{align*} 
Then by the previous lemma, there exists a random ball in $V$ which absorbs the bounded sets of $H$. Since $V$ is compactly embedded in $H$, there exists a compact set $K\subset H$ such that, for all bounded set $B\subset H$ there exists $\bar{t}\le -1$ such that $\varphi B\subset K$ $\P$ a.s.. Defining $K(\omega) :=\{u\in H: |u|\le r_2(\omega)\}$, we have proved the existence of a compact absorbing set. Then by Theorem \ref{compactattract}, there exists random attractor to (\ref{asnse1rds}).
The existence of an invariant Markov measure is a direct consequence of Corollary \ref{fpm}, provided we can show that the one-point motions associated with the flow $\varphi(t,\omega)$ define a family of Markov processes. The proof of this is analogous to the proof of Markov property of solutions to the (\ref{asnse1rds}) in \cite{snse1}. Nevertheless we repeat here as well.
Let $\varphi_{s,t}$ be defined as in earlier section. Let $\mathcal{F}_{s,t}$ be the $\sigma$-algebra generated by $L(r)-L(s)$ for all $r\in[s,t]$, and let $\mathcal{F}_t=\mathcal{F}_{0,t}.$ Define the operators $P_t$ in the space of bounded measurable function over $H$ as $(P_t f)(u_0)=\E f(\varphi(t)u_0).$ 
To prove $\varphi(t,\omega)$ defines a family of Markov processes. It suffices to prove 
\begin{align*}
    \E[f(\varphi(t+s)x)|\mathcal{F}_t]=P_{s}(f)(\varphi(t)x),
\end{align*}   
for all $0\le s\le t$ and all bounded continuous functions $f$ over $H$, which implies that $\varphi(t+s)x$ is a Markov process with transition semigroup $P_t$.  By uniqueness, the following holds
 \begin{align*}     \varphi(t+s,\omega)x=\varphi(s,\omega)\varphi(t,\omega)x 
 \end{align*}  
over $[t,\infty]$ with $\mathcal{F}_t$ measurable initial condition $\varphi(t,t)\delta=\delta.$
It suffices to prove 
\begin{align}\label{psfe}
    \E [f(\varphi_{t,t+s}\delta)|\mathcal{F}_t]=P_s(f)(\delta)
\end{align}
 for for every 
 $H$ integrable, $\mathcal{F}_t$ random variable $\delta$.

Note, (\ref{psfe}) not only holds for every $f\in C_b(H)$, but also holds for $\varphi=1_{\Gamma}$, where $\Gamma$ is an arbitrary Borel set of $H$ and consequently  for all $\varphi\in B_b(H)$. Without loss of generality, we assume $\varphi\in C_b(H).$ We know that, if $\delta=\delta_i$ $\P$ a.s., then $random variable\, \varphi(t,t+s)\delta_i$ is independent to $\mathcal{F}_t$, since $\varphi(t,t+s)\delta_i$  is $\mathcal{F}_{t,t+s}$ measurable. Hence,
\begin{align*}
    \E(f(\varphi(t,t+s)\delta_i)|\mathcal{F}_t)=\E f(\varphi(t,t+s)\delta_i)=P_{t,t+s}f(\delta_i)=P_s f(\delta_i),\quad \P\quad\text{a.s.}\quad .
\end{align*}   
Since the coefficient of the equation for $\varphi(t,t+s)$ are independent, one can see that the $H$ random variable $\varphi_{t,t+s}$ and $\varphi_s x$ have the same law. If $\delta$ has the form
\begin{align}\label{simpleeta}  \delta=\sum^N_{i=1}\delta^{i}1_{\Gamma^i},
\end{align}
where $\delta^{(i)}\in H$ and  $\Gamma^{(i)}\subset\mathcal{F}_t$ is a partition of $\Omega$, $\delta_i$ are elements of $H$. Then
\begin{align*}
    \varphi(t,t+s)\delta_i=\sum^N_{i=1}\varphi(t,t+s,\delta_i)1_{\Gamma_i},\quad\P,\quad\text{a.s.}\quad .
\end{align*}
Hence,
\begin{align*}
    \E(f(\varphi(t,t+s)\delta)|\mathcal{F}_{t})=\sum^N_{i=1}\E(f(\varphi(t,t+s)\delta_i)1_{\Gamma_i}|\mathcal{F}_t)\quad\P\quad\text{a.s.}\quad .
\end{align*}
Take into account the random variable $u(t,t+s)\delta_i$ independent to $\mathcal{F}_t$ and $1_{\Gamma_i}$ are $\mathcal{F}_t$ measurable, $i=1,\cdots, l$, one deduces that 
\begin{align*}
    \E[f(\varphi(t,t+s,\delta))|\mathcal{F}_t]&=\sum^N_{i=1} P_s f(\delta_i)1_{\Gamma_i}=P_s f(\delta),\quad\P\quad\text{a.s.}
\end{align*}
and so (\ref{psfe}) is proved. For a general $\delta$ there exists a sequence of $\delta_n$ for which (\ref{psfe}) holds converges to $\delta$ in $L^2(\Omega;H)$ a.s., that is,
\begin{align*}
    \E|\delta-\delta_n|^2\to 0.
\end{align*}
By continuity of $f$ one can pass in the identity (\ref{psfe}), with $\delta$ replaced with $\delta_n$, to the limit and (\ref{psfe}) holds if $\E|\delta|^2<\infty$. So $\varphi(t,\omega)$ defines a family of Markov processes.

The proof of existence of Markov measure is completed.
\end{proof}
\begin{rmk}
    Although the same results hold in $\beta$-stable L\'evy case as in the Gaussian case (see \cite{MR1305587}), there is some difference between dealing with Brownian motion and L\'evy motions. First, we need to consider c\`adl\`ag function in the Skorohod metric, which are different from the continuous case in the metric under the compact-open topology. Second, one has to consider solutions in the sense of Carth\'eodory and the right-hand derivatives.  
\end{rmk} 


Let $u(t,x)$ be the unique solution to problem (\ref{asnse1rds}). Let us recall from \cite{snse1} that such a unique solution exist for each $x\in H.$ Let us define the transition operator $P_t$ by a standard formula. For $f\in C_b(H),$ put 
\begin{align*}
    (P_t f)(x)=\E f(\varphi(t,x)),\quad t\ge 0, x\in X.
\end{align*}
In view of Proposition \ref{3.8}, $(P_t, t\ge 0)$ is a family of Feller operators, i.e. $P_t: C_b(H)\to C_b(H)$ and, for any $f\in C_b(H)$ and $x\in H$, $P_t f(x)\to f(x)$ as $t\downarrow 0.$
Moreover, following the identical lines of the proof of Theorem \ref{t5attr} in last subsection, one can prove that $\varphi$ is a Markov RDS. Invoke corollary \ref{fpm}, we deduce the existence of Feller invariant measure for our stochastic Navier-Stokes equations (\ref{asnse1rds}).

\begin{cor}\label{imrdsbrze}
    There exists an Feller Markov Invariant Measure for the SNSE (\ref{asnse1rds})
\end{cor}

\section{Acknowledgments}
This work is taken out of the author's PhD thesis (awarded in April 2018). The author would like to express her gratitude to her PhD advisor Prof. Ben Goldys for his helpful suggestions to this work. The author is grateful to Prof. Zdzislaw Brzezniak for his insightful comments. 
\bibliographystyle{apalike}
\bibliography{l_sthesis}             
\end{document}